# Heavy tail properties of stationary solutions of multidimensional stochastic recursions

## Yves Guivarc'h[1]

*Université de Rennes 1 (France)*

**Abstract:** We consider the following recurrence relation with random i.i.d. coefficients $(a_n, b_n)$:

$$x_{n+1} = a_{n+1}x_n + b_{n+1}$$

where $a_n \in GL(d, \mathbb{R}), b_n \in \mathbb{R}^d$. Under natural conditions on $(a_n, b_n)$ this equation has a unique stationary solution, and its support is non-compact. We show that, in general, its law has a heavy tail behavior and we study the corresponding directions. This provides a natural construction of laws with heavy tails in great generality. Our main result extends to the general case the results previously obtained by H. Kesten in [16] under positivity or density assumptions, and the results recently developed in [17] in a special framework.

## 1. Notation and problem

### 1.1. General notation

We consider the $d$-dimensional vector space $V = \mathbb{R}^d$, endowed with the scalar product

$$\langle x, y \rangle = \sum_1^d x_i y_i$$

and the corresponding norm

$$\|x\| = \left( \sum_1^d |x_i|^2 \right)^{1/2}.$$

Let $G = GL(V)$ be the general linear group, and $H = Aff(V)$ the affine group. An element $h \in H$ can be written in the form

$$h(x) = ax + b \quad (x \in V)$$

where $a \in G$, $b \in V$. In reduced form we write $h = (b, a)$ and we observe that the projection map $(b, a) \to a$ is a homomorphism from $H$ to $G$. We consider also the projection $(b, a) \to b$ of $H$ on $V$ and we observe that $H$ can be written as a semi-direct product $H = G \ltimes V$ where $V$ denotes also the translation group of $V$.

For a locally compact second countable (l.c.s.c.) space $X$, we denote by $M^1(X)$ the convex set of probability measures on $X$ and we endow $M^1(X)$ with the weak topology. For a l.c.s.c. semi-group $L$ and a probability measure $\rho \in M^1(L)$, we denote by $T_\rho$ the closed semi-group generated by *supp* $\rho$, the support of $\rho$. In particular let $\eta \in M^1(H)$ and denote by $\mu$ (resp. $\overline{\eta}$) its projection on $G$ (resp. $V$). Then $T_\eta$ (resp. $T_\mu$) will denote the corresponding semi-group in $H$ (resp. $G$). Let $L$ (resp. $X$) be a l.c.s.c. semi-group (resp. l.c.s.c. space) and assume $X$ is a $L$-space. Then for any $\rho \in M^1(L)$ and $\sigma \in M^1(X)$, we define the convolution $\rho * \sigma \in M^1(X)$ by:

$$\rho * \sigma = \int \delta_{hx} d\rho(h) \, d\sigma(x) = \int h_*(\sigma) \, d\rho(h)$$

where $h_*(\sigma) \in M^1(X)$ is the push-forward of $\sigma$ under the map $h \in L$. We are interested in the special cases $L = H$ or $G$, $X = V$ or $V \setminus \{0\}$ and the actions of $H$, $G$ are the natural ones. We need to consider a class of Radon measures $\lambda$ on $V \setminus \{0\}$ which are homogeneous under dilations. Such a measure can be written in the form $\lambda = \rho \otimes \ell^s$ where $\rho$ is a measure on $\mathbb{S}^{d-1}$ and $\ell^s(dt) = \frac{dt}{t^{s+1}}$ is the natural s-homogeneous measure on $\mathbb{R}_+ = \{t \in \mathbb{R} \; ; \; t > 0\}$.

### 1.2.   *The recursion*

We consider the product space $\Omega = H^{\mathbb{N}}$ endowed with the product measure $\mathbb{P} = \eta^{\otimes \mathbb{N}}$ and the shift map $\theta$. The corresponding expectation is denoted by $\mathbb{E}$. We are interested in the stationary solutions of the equation:

$$(1.0) \qquad\qquad x_{n+1} = a_{n+1} x_n + b_{n+1}$$

where $a_n \in GL(d, \mathbb{R}), b_n \in \mathbb{R}^d$, i.e., sequences $(x_n)_{n \in \mathbb{N}}$ such that the law of $(x_n)_{n \in \mathbb{N}}$ is shift-invariant. As is well known from the theory of Markov chains, this reduces to finding probability measures $\nu \in M^1(V)$ (the law of $x_0$) which are $\eta$-stationary, i.e.:

$$(1.1) \qquad\qquad \nu = \eta \, * \, \nu$$

Here, we are interested in "the shape at infinity" of $\nu$, if such a $\nu$ exists. A sufficient condition for the existence of $\nu$ is the convergence of the random series $z = \sum_0^\infty a_1 \cdots a_k b_{k+1}$. Then $z$ satisfies the functional equation:

$$(1.2) \qquad\qquad z = a_1 \, (z \circ \theta) + b_1$$

and $\nu$ is the law of $z$.

In particular, the above series converges if $\mathbb{E}(\log^+ \|a\| + \log^+ \|b\|) < +\infty$ and $\alpha = \lim_k \frac{1}{k} \mathbb{E}(\log \|a_0 \cdots a_k\|) < 0$. We denote by

$$\mathbb{S}^{d-1} = \{u \in V \; : \; \|u\| = 1\}$$

the unit sphere in $\mathbb{R}^d$, and for $u \in \mathbb{S}^{d-1}, t \geq 0$, we write

$$H_t^u = \{v \in V \; : \; \langle u, v \rangle \geq t\}.$$

Denote the Radon transform of $\nu$ by $\varphi_u(t)$, and the asymptotic behavior of $\varphi_u(t) = \nu(H_t^u)$ (as $t \to +\infty$) is the so-called tail of $\nu$ in the direction $u$. In order to describe $\nu$, we compactify the Euclidean space $V$ by adding the sphere at infinity $\mathbb{S}_\infty^{d-1}$: $\overline{V} = V \cup \mathbb{S}_\infty^{d-1}$. A sequence $v_n \in V$ converges to $u \in \mathbb{S}_\infty^{d-1}$ if

$$\lim_n \|v_n\| \to \infty \quad \text{and} \quad \lim_n \frac{v_n}{\|v_n\|} = u,$$

where $u \in \mathbb{S}^{d-1}$ is naturally identified with the point at infinity in $\mathbb{S}_\infty^{d-1}$ with the same direction. The space $V$ is endowed with the metric

$$d(v, v') = \|v - v'\|,$$

and for $Y \subset V, x \in V$ we denote

$$d(x, Y) = \inf_{y \in Y} \|x - y\|.$$

In a sense, the problem reduces to finding the asymptotics of $\varphi_u(t)$ for large $t$. For the case of positive matrices $a_n$ and positive vectors $b_n$, it was shown by H. Kesten in [16] that, under some non-degeneracy conditions, there exist $\chi > 0$ and $c > 0$ such that

$$\lim_{t \to +\infty} t^\chi \nu(B_t') = c > 0,$$

where $B_t' = \{x \in V : \|x\| \geq t\}$. Furthermore, there exists a function $c(u)$ on $\mathbb{S}^{d-1}$ which is not identically zero, such that

$$\forall \, u \in \mathbb{S}^{d-1}, \, \lim_{t \to +\infty} t^\chi \nu(H_t^u) = c(u).$$

Under some assumption on $\eta$ to be detailed below, our main result implies in particular the following:

**Theorem 1.1.** *There exist $\chi > 0$ and a non-zero positive measure $\nu_\infty$ on $\mathbb{S}^{d-1}$, such that*

$$\lim_{t \to +\infty} t^\chi (\delta_{t^{-1}} * \nu)(\varphi) = (\nu_\infty \otimes \ell^\chi)(\varphi)$$

*for every Borel function $\varphi$, where the set of discontinuities of $\varphi$ has $\nu_\infty \otimes \ell^\chi$ measure $0$, and for some $\varepsilon > 0$,*

$$\sup_{v \neq 0} \|v\|^{-\chi} \, |\log \|v\||^{1+\varepsilon} \, |\varphi(v)| < +\infty.$$

Under irreducibility conditions, we will describe in Section 6 the properties of $\nu_\infty$, depending on the geometry of $T_\eta$ and $T_\nu$. A subtle point is the calculation of the support of $\nu_\infty$. For $d = 1$, this reduces to the study of the positivity of $c_+$ and $c_-$ (see [9] and Theorem 3.1). Under special hypothesis, the calculation of the support of $\nu_\infty$ is done in [17], which gives $supp \, \nu_\infty = \mathbb{S}^{d-1}$ (see Corollary 6.4). In particular, Theorem 6.2 and Theorem 3.1 below can be considered as extensions of the results in [9], [16], [17]. The method of proof that we will use is strongly inspired by [16]. It can be roughly described as "linearization at infinity" of the stationarity equation $\nu = \eta * \nu$, (see Section 4). More precisely, this equation will be replaced by

(1.3) $$\lambda = \mu * \lambda .$$

Thus the problem reduces to comparing $\nu$ and $\lambda$ at infinity. In a sense to be explained below, (1.3) is the homogeneous equation associated with the inhomogeneous linear equation (1.1). We will be able to express $\nu$ in terms of the Green kernel of equation (1.3), i.e., the potential kernel of $\mu$ (see Section 7). It is known that equation (1.3) has only one or two extremal solutions which are relevant to the problem. In the proofs, we will rely strongly on the analytic tools of [18], which have also been essential in [17]. We will also use the dynamical aspects of linear group actions on real vector spaces (see the recent surveys [13], [14]). For information on products of random matrices, we refer to [2], [11], [12] and [13]. For an account of our main result in a special case, see [6].

The present exposition is on results which extend the results in the joint work [6]. It is also an improved version of part of the unpublished work [18].

## 2. The stationary measure and its support

### 2.1.  *Notation*

For any $g \in G$, we denote

$$r(g) = \lim_n \|g^n\|^{1/n}.$$

We say that $g \in G$, or $h = (b, a) \in H$, is quasi-expanding (resp. contracting) if $r(g) > 1$ (resp. $r(g) < 1$). For a semi-group $T \subset G$ or $T \subset H$, we will denote by $T^e \subset T$ (resp. $T^c \subset T$) the subset of its quasi-expanding (resp. contracting) elements. We observe that, if $h = (b, a)$ is contracting, then $h$ has a unique (attractive) fixed point $h^a \in V$:

$$h^a = (I - a)^{-1} b = \sum_0^\infty a^k b \ .$$

For a semi-group $T \subset H$ we denote:

$$\Lambda^a(T) = \text{Closure}\{h^a : \ h \in T^c\} \subset V \ .$$

In what follows, a measure $\eta \in M^1(H)$ will be given and the following moment conditions will be assumed:

$$(M) \qquad \int \log^+ \|a\| d\mu(a) < +\infty \quad \text{and} \quad \int \log^+ \|b\| d\overline{\eta}(b) < +\infty \ .$$

We will consider the largest Lyapunov exponent $\alpha$ of the random product $a_0 \cdots a_n$:

$$\alpha = \lim_n \frac{1}{n} \int \log \|a\| d\mu^n(a) \ ,$$

where $\mu^n$ is the $n^{th}$ convolution power of $\mu$ and the limit exists by subadditivity. We observe that, if $\alpha < 0$, then there exists $a \in T_\mu$ with $\|a\| < 1$, and hence $r(a) < 1$. It follows that $T_\mu^c \neq 0$ and $\Lambda^a(T_\eta) \neq \phi$.

### 2.2.  *Existence and uniqueness of $\nu$*

As a preliminary result we state the following, the first part of which is well known [3].

**Proposition 2.1.** *Assume that $\eta \in M^1(H)$ satisfies condition*

$$(M) \qquad \qquad \mathbb{E}(\log^+ (\|a\|)) < +\infty \quad \text{and} \quad \mathbb{E}(\log^+ \|b\|) < +\infty \ ,$$

*and* $\alpha = \lim_n \frac{1}{n} \mathbb{E}(\log \|a_0 \cdots a_n\|) < 0$. *Assume also that* $T_\eta$ *has no fixed point. Then there exists a unique stationary solution of* $\nu = \eta * \nu$.

*The probability measure* $\nu$ *is non-atomic, the series* $\sum_0^\infty a_1 \cdots a_k b_{k+1}$ *converges* $\mathbb{P} - a.s.$ *with* $\nu$ *as its limiting law. The support of* $\nu$ *is equal to* $\Lambda^a(T_\eta)$ *and* $\Lambda^a(T_\eta)$ *is the unique* $T_\eta$*-minimal set in* $V$.

*For every initial point* $x_0 = x$, *the random vector* $x_n$ *defined by*

$$x_{n+1} = a_{n+1} x_n + b_{n+1}$$

*converges in law to* $\nu$ *and we have the* $\mathbb{P} - a.s.$ *convergence*

$$\lim_n d[x_n, \Lambda^a(T_\eta)] = 0 .$$

*Furthermore, if* $T_\mu^e \neq \phi$, *then* $\Lambda^a(T_\eta)$ *is non-compact.*

## 2.3. Some examples

As an illustration, we consider some special examples.

We first consider the special case $d = 1$, $b_k = 1$ and $a_k > 0$, i.e. the recursion

$$x_{n+1} = a_{n+1} x_n + 1 .$$

If $\mathbb{E}(\log a) < 0$, then its stationary solution has the same law as

$$z = 1 + \sum_1^\infty a_1 \cdots a_k .$$

We take $\mu$ of the form $\mu = \frac{1}{2}(\delta_u + \delta_{u'})$ with $\alpha = \mathbb{E}(\log a) = \frac{1}{2}(\log u + \log u') < 0$.

(a) Choose $u = \frac{1}{2}, u' = \frac{1}{3}$, then $\alpha = -\frac{1}{2}\log 6 < 0$, $T_\mu^e = \phi$ and $supp\ \nu$ is a Cantor subset of the interval $[\frac{3}{2}, 2]$. Hence, the tail of $\nu$ vanishes.

(b) Choose $u = 1/3$, $u' = 2$, then $\alpha = -\frac{1}{2}\log(3/2) < 0$ and $2 \in T_\mu^e \neq \phi$. It follows from [16] (see also [5], [7], [9], [10]) that

$$supp\ \nu = [\frac{3}{2}, \infty[ \quad \text{and} \quad \lim_{t \to +\infty} t^\chi \nu(t, \infty) = c > 0 ,$$

where $\chi > 0$ is defined by $\mathbb{E}(a^\chi) = 1$. Such a $\chi$ exists since the log-convex function $k(s) = \mathbb{E}(|a|^s)$ satisfies

$$\alpha = k'(0) < 0 \quad \text{and} \quad \lim_{s \to +\infty} k(s) = +\infty.$$

An example in two dimension was proposed by H. Kesten in [16]. A simplified version is the following

$$\eta = p\delta_h + p'\delta_{h'} , \qquad p, p' > 0 ,$$

$$h = \rho \begin{pmatrix} \cos\theta & -\sin\theta \\ \sin\theta & \cos\theta \end{pmatrix} , \quad h' = \begin{pmatrix} b, \begin{pmatrix} \lambda & 0 \\ 0 & \lambda' \end{pmatrix} \end{pmatrix} ,$$

$$\theta/\pi \notin \mathbb{Q}, \quad b \neq 0, \quad 0 < \rho < 1, \quad 0 < \lambda' < 1 < \lambda.$$

Since $\lambda > 1$, we have $h' \in T_\mu^e \neq \phi$. Also, if $\rho$ is sufficiently small, then we have $\alpha < 0$.

It will follow from our main result (Section 6) that for any $u$ with $\|u\| = 1$, there exist $\chi > 0$ and $c(u) > 0$, such that

$$\lim_{t \to +\infty} t^\chi \nu(H_t^u) = c(u) > 0 .$$

### 3. The case of the line ($d = 1$)

For the multidimensional case ($d > 1$), we will impose geometric conditions on $T_\mu$, which imply non-arithmeticity properties. The case $d = 1$ will not be covered by these general assumptions. Furthermore, for $d = 1$, we will remove the condition $\det(a) \neq 0$, which means that we replace the group $H$ by the affine semi-group $H_1$ of $V$:

$$H_1 = \{g = (b, a) \; ; \; b \in \mathbb{R} \; , \; a \in \mathbb{R}\}.$$

However the final results will be similar. As a comparison and introduction to the general case, we first give the result for $d = 1$. We need to consider the Mellin transform of $\mu$:

$$k(s) = \int |a|^s d\mu(a) \; .$$

We suppose that $k(s)$ is defined for some $s > 0$, i.e.,

$$s_\infty = \sup\{s \geq 0; k(s) < \infty\} > 0.$$

Our main condition on $\mu$, which is responsible for the tail behavior of $\nu$, will be the following contraction-expansion condition:

$$(C\text{--}E) \qquad \alpha = k'(0) < 0, \qquad s_\infty > 0, \qquad \lim_{s \to s_\infty} k(s) \geq 1.$$

It is satisfied if $s_\infty = +\infty$, $\alpha < 0$ and $T_\mu \not\subset [-1, 1]$. Using condition $(C - E)$, we can define $\chi > 0$ by $k(\chi) = 1$.

   We also need to consider the corresponding probability measure $\mu^\chi(da) = |a|^\chi \mu(da)$. We denote by $(C)$ the following set of conditions:

$$(M) \qquad\qquad \int (\log^+ |a| + \log^+ |b|) d\eta(h) < +\infty \; ,$$

$$(C\text{--}E) \qquad\qquad \alpha < 0, \qquad s_\infty > 0, \qquad \lim_{s \to s_\infty} k(s) \geq 1,$$

$$(M_\chi) \qquad\qquad \int [\, |a|^\chi \log^+ |a| + |b|^\chi] d\eta(h) < +\infty \; ,$$

$$(N\text{--}F) \qquad\qquad\qquad T_\eta \text{ has no fixed point.}$$

**Theorem 3.1** ($H_1 = "ax + b"$). *Assume that $\eta \in M^1(H_1)$ satisfies hypothesis $(C)$, and for any $\rho > 0$,*

$$(N\text{--}A) \qquad\qquad supp \; \mu \not\subset \{\pm\rho^n : n \in \mathbb{Z}\} \cup \{0\},$$

*then $\nu$ is diffuse and there are only 3 cases.*

*Case I: supp $\mu \not\subset \mathbb{R}_+ \cup \{0\}$,*
   *then supp $\nu = \mathbb{R}$ and there exists $c > 0$ such that*

$$\lim_{t \to +\infty} t^\chi \nu(t, \infty) = \lim_{t \to -\infty} |t|^\chi \nu(-\infty, t) = c.$$

*Case II: supp $\mu \subset \mathbb{R}_+ \cup \{0\}$.*

*Case II 1:* $+\infty \in \overline{\Lambda^a(T_\eta)}$ *and* $-\infty \in \overline{\Lambda^a(T_\eta)}$,

then *supp* $\nu = \mathbb{R}$, *and there exist* $c_+ > 0$, $c_- > 0$ *such that*

$$\lim_{t \to +\infty} t^\chi \nu(t, \infty) = c_+ \quad and \quad \lim_{t \to -\infty} |t|^\chi \nu(-\infty, t) = c_-.$$

*Case II 2:* $+\infty \in \overline{\Lambda^a(T_\eta)}$ *but* $-\infty \notin \overline{\Lambda^a(T_\eta)}$,

then *there exist* $c > 0$ *and* $m \in \mathbb{R}$ *such that*

$$supp \; \nu = [m, +\infty) \quad and \quad \lim_{t \to +\infty} t^\chi \nu(t, \infty) = c.$$

**Remark.** For previous work on the case $d = 1$, see [9], where it is proved that $c_+ + c_- > 0$ and expressions for $c_+$ and $c_-$ are given in terms of $\nu$. Our results on *supp* $\nu$ and the tail of $\nu$ in case II 1 are new. For an application to tail estimates in a different context, see [10].

## 4. The linearization procedure ($d > 1$)

Here we develop a heuristic approach which suggests that, at infinity, the solution $\nu$ of $\nu = \eta * \nu$ should be compared to a "stable solution at infinity" of the linear homogeneous equation $\lambda = \mu * \lambda$, where $\lambda$ is a Radon measure on $V \setminus \{0\}$. We think of $\nu$ as a perturbation of the trivial solution $\delta_0$ of the unperturbed equation $\delta_0 = \mu * \delta_0$.

### *4.1. Homogeneity at infinity*

For $s \geq 0$, we denote

$$
\begin{aligned}
k(s) &= \lim_n \left[ \int \|g\|^s d\mu^n(g) \right]^{1/n}, \\
s_\infty &= \sup\{s \geq 0 : k(s) < +\infty\},
\end{aligned}
$$

and we assume

$$s_\infty > 0, \quad \lim_{s \to s_\infty} k(s) \geq 1, \quad \alpha = \lim_n \frac{1}{n} \int \log \|g\| d\mu^n(g) < 0.$$

We then define $\chi > 0$ by $k(\chi) = 1$. It follows that if $s \in (0, \chi)$, then $k(s) < 1$. We estimate the $s$-th moment of $z = \sum_0^\infty a_1 \cdots a_k b_{k+1}$:

$$\mathbb{E}(\|z\|^s) \leq \sup(1, s) \sum_0^\infty \mathbb{E}(\|a_1 \cdots a_k b_{k+1}\|^s).$$

Since $\lim_n \mathbb{E}(\|a_1 \cdots a_k b_{k+1}\|^s)^{1/n} = k(s) < 1$ for $s \in (0, \chi)$, we conclude that $\mathbb{E}(\|z\|^s) < +\infty$. This calculation is not valid for $s = \chi$. This suggests that $\mathbb{E}(\|z\|^\chi) = +\infty$, and that if $\eta$ is sufficiently non-degenerate, then $\nu$ will have the following homogeneity property at infinity:

$$\lim_{t \to +\infty} t^\chi (\delta_{t^{-1}} * \nu) = \lambda \; ,$$

where $\lambda$ is $\chi$-homogeneous under dilations, i.e., $\lambda = \nu_\infty \otimes \ell^\chi$ for some non-zero measure $\nu_\infty$ on $\mathbb{S}^{d-1}$. Furthermore, if the convergence holds for bounded continuous functions, then $\nu_\infty$ will satisfy:

$$\lim_{t \to +\infty} t^\chi \nu_{B'_t} = \nu_\infty \quad \text{and} \quad \lim_{t \to +\infty} t^\chi \nu(H^u_t) = \nu_\infty(\overline{H}^u_0 \cap \mathbb{S}^{d-1}_\infty),$$

where $\nu_{B'_t}$ is the restriction of $\nu$ to $B'_t = \{v \in V : \|v\| \geq t\}$, $\nu_\infty$ is a positive measure on $\mathbb{S}^{d-1}_\infty$ identified with $\mathbb{S}^{d-1}$, and $\overline{H}^u_0$ is the closure of $H^u_0$ in $\overline{V}$.

### *4.2.   Derivation of the linearized equation*

We restrict to the case where $\eta$ is a product measure on $H$. In this case,

$$\eta = \overline{\eta} * \mu \ , \qquad \nu = \overline{\eta} * (\mu * \nu).$$

Writing $\nu_t = t^\chi(\delta_{t^{-1}} * \nu)$ and $\overline{\eta}^t = \delta_{t^{-1}} * \overline{\eta} * \delta_t$, we get $\nu_t = \overline{\eta}^t * (\mu * \nu_t)$.

Since $\lim\limits_{t \to +\infty} \overline{\eta}^t = \delta_0$ and $\lim\limits_{t \to +\infty} \nu_t = \lambda$, we obtain

$$\lambda = \mu * \lambda, \qquad \lambda = \nu_\infty \otimes \ell^\chi.$$

We introduce the natural actions of $g \in G$ on $\mathbb{S}^{d-1}$ and $\mathbb{P}^{d-1} = \mathbb{S}^{d-1}/\{\pm Id\}$. On $\mathbb{S}^{d-1}$ we have $g.u = \frac{gu}{\|gu\|}$, $(u \in \mathbb{S}^{d-1})$. The projective action of $g$ on $x \in \mathbb{P}^{d-1}$ will also be denoted $g.x$.

Then the above equation reduces to

$$\nu_\infty = \int \|gx\|^\chi \delta_{g.x} d\nu_\infty(x) d\mu(g),$$

and a similar equation on $\mathbb{P}^{d-1}$ holds for the projection of $\nu_\infty$. Equations of this type were considered by H. Furstenberg (see [8], [11], [13]) in the context of harmonic measures for random walks.

We introduce the representations $r^s$ on $\mathbb{S}^{d-1}$, $\mathbb{P}^{d-1}$:

$$r^s(g)(\delta_x) = \|gx\|^s \delta_{g.x} \ .$$

In particular, the above integral equation can be written as

$$r^\chi(\mu)(\nu_\infty) = \nu_\infty \ .$$

## 5. **Limit sets on $\mathbb{P}^{d-1}$, $\mathbb{S}^{d-1}$**

### *5.1.   Notation*

We recall briefly some definitions and results of [11], [13].

**Definition 5.1.** A semi-group $T \subset GL(V)$ is said to be strongly irreducible if there does not exist a finite union of proper subspaces of $V$ which is $T$-invariant.

**Definition 5.2.** An element $g \in GL(V)$ is said to be proximal if there exists a unique eigenvalue $\lambda_g$ of $g$, such that

$$|\lambda_g| = r(g) = \lim_n \|g^n\|^{1/n} \ .$$

This means that we can write

$$V = \mathbb{R}v_g \oplus V_g^< \ ,$$

where $gv_g = \lambda_g v_g$, $V_g^<$ is $g$-invariant, and the spectral radius of $g$ in $V_g^<$ is strictly less than $|\lambda_g|$. We denote by $g^a \in \mathbb{P}^{d-1}$ the point corresponding to $v_g \in V$. We observe that $g^a$ is attractive:

$$\forall \ x \notin V_g^<, \quad \lim_n g^n.x = g^a \ .$$

**Definition 5.3.** We say that $T \subset G$ satisfies condition i.p. if $T$ is strongly irreducible and contains a proximal element. We denote by $T^{prox}$ the subset of proximal elements.

If for example, the Zariski closure of $T$ contains $SL(V)$, then it is known that $T$ satisfies i.p. (see [12], [13]). Furthermore, condition i.p. for $T$ is equivalent to condition i-p for the Zariski closure of $T$ [13]. This Zariski closure is a Lie group with a finite number of components with a special structure described in [14], Lemma 2.7.

### 5.2. The dynamics of $T$ on $\mathbb{P}^{d-1}, \mathbb{S}^{d-1}$

Condition i.p. ensures that the dynamics of $T$ on $\mathbb{P}^{d-1}$ and $\mathbb{S}^{d-1}$ can be described in a simple way.

**Definition 5.4.** Assume $T \subset GL(V)$ satisfies condition i.p. We denote

$$\Lambda(T) = \text{Closure}\{g^a \in \mathbb{P}^{d-1} : g \in T^{prox}\} \ ,$$
$$\Lambda_1(T) = \text{Closure}\{v_g \in \mathbb{S}^{d-1} : g \in T^{prox}\} \ .$$

**Proposition 5.5.** *Assume that* $T \subset GL(V)$ *satisfies condition i.p. of Definition 5.3. Then* $\Lambda(T)$ *is the unique* $T$-*minimal set of* $\mathbb{P}^{d-1}$. *The action of* $T$ *on* $\mathbb{S}^{d-1}$ *has either one or two minimal sets, whose union is* $\Lambda_1(T)$.
*Case I: $T$ does not preserve a convex cone in $V$,*
  *then $\Lambda_1(T)$ is the unique $T$ minimal set.*
*Case II: $T$ preserves a convex cone $\mathcal{C} \subset V$,*
  *then the action of $T$ on $\mathbb{S}^{d-1}$ has two minimal sets $\Lambda_1^+(T), \Lambda_1^-(T)$ with*

$$\Lambda_1^+(T) \subset \mathbb{S}^{d-1} \cap \mathcal{C} \quad , \quad \Lambda_1^-(T) = -\Lambda_1^+(T).$$

The existence of a convex cone preserved by $T$ is not related to the fact that $T$ is a semigroup and not a group. Examples of Zariski dense groups preserving a convex cone exist in abondance (See [1]). For the action of $\mu$ on measures we have the following

**Proposition 5.6.** *Assume $T = T_\mu$ satisfies condition i.p., $s_\infty > 0$, $\lim\limits_{s \to s_\infty} k(s) \geq 1$. Let $\chi > 0$ be defined by $k(\chi) = 1$. Then the equation*

$$r^\chi(\mu)(\rho) = \rho, \quad \rho \in M^1(\mathbb{S}^{d-1})$$

*has one or two extremal solutions.*
  *In case I as above, $\rho = \nu_1$.*
  *In case II, there are two extremal solutions $\nu_1^+, \nu_1^-$ with*

$$supp \ (\nu_1^+) = \Lambda_1^+(T_\mu), \quad supp \ (\nu_1^-) = \Lambda_1^-(T_\mu) \ ,$$

*and $\nu_1^-, \nu_1^+$ are symmetric with respect to each other.*

An important consequence of i.p. which guarantees the $\chi$-homogeneity of $\lambda$ at infinity is given by the following (see [11], [14]). Together with Proposition 5.6, it is one of the main algebraic facts which play a role in Theorem 6.2 below.

**Proposition 5.7** ($d > 1$). *Assume that $T \subset GL(V)$ satisfies i.p. Then the subgroup of $\mathbb{R}$ generated by the spectrum of $T$,*

$$\Sigma(T) = \{\log |r(g)| \; ; \; g \in T^{prox}\},$$

*is dense in $\mathbb{R}$.*

## 6. The main theorem ($d > 1$)

### 6.1. Main theorem

As in Section 3, we will use the following set of hypothesis $(C)$:

$(M)$ $$\int [\log^+ \|a\| + \log^+ \|b\|] d\eta(h) < +\infty \; ,$$

$(C–E)$ $$s_\infty > 0, \quad \alpha < 0, \quad \lim_{s \to s_\infty} k(s) \geq 1 \; ,$$

$(M_\chi)$ $$\int \|a\|^\chi \log^+ \|a\| d\mu(a) < +\infty, \quad \int \|b\|^\chi d\overline{\eta}(b) < +\infty \; ,$$

$(N–F)$ $$T_\eta \text{ has no fixed point} \; .$$

We observe that $k(s)$ is a log-convex function and is finite at $s$ if $\int \|g\|^s d\mu(g) < +\infty$. Also, if $\alpha < 0$ and $\lim_{s \to s_\infty} k(s) \geq 1$, then there exists $\chi > 0$ with $k(\chi) = 1$. If $\alpha < 0$, $T_\mu$ contains a quasi-expanding element and $s_\infty = +\infty$, then $\lim_{s \to +\infty} k(s) = +\infty$. A detailed study of $k(s)$ under condition i.p. can be found in [13], where $k(s)$ is shown to be analytic. Condition i-p will also be used in Theorem 6.2 below.

The above conditions are inspired by [16], where the case of positive matrices is considered. There, it is not assumed that $\det(a) \neq 0$, and reducibility is allowed. However, hypothesis $(C)$ and existence of a proximal element are implicitly assumed. On the other hand, for $d > 1$, the non-arithmeticity condition of [16] does not appear explicitly here, although it is valid (compare with Theorem 3.1).

We will identify the sphere at infinity, $\mathbb{S}^{d-1}_\infty$, with the unit sphere $\mathbb{S}^{d-1}$. In particular we consider $\Lambda_1^+(T_\mu)$, $\Lambda_1^-(T_\mu)$ and $\Lambda_1(T_\mu)$ as subsets of $\mathbb{S}^{d-1}_\infty$. The convex envelope of a closed subset $Y \subset \mathbb{R}^d$ will be denoted by $Co(Y)$. If $Y \subset \mathbb{S}^{d-1}$, we also denote by $Co(Y)$ the intersection of $\mathbb{S}^{d-1}$ with the convex envelope of the cone generated by $Y$.

We will need a concept of direct Riemann integrability as in [16], [17].

**Definition 6.1.** Let $X$ be a compact metric space, $\rho$ a probability measure on $X$, and $\varphi$ a Borel function on $X \times \mathbb{R}_+$. We say that $\varphi$ is $\rho$-directly Riemann integrable if

(a) $$\sum_{k=-\infty}^{+\infty} \sup\{|\varphi(x,t)| : (x,t) \in X \times [2^k, 2^{k+1}]\} < +\infty \; ,$$

(b) The set of discontinuities of $\varphi$ has $\rho \otimes \ell$ measure 0.

Using polar coordinates, we can write $V \setminus \{0\} = \mathbb{S}^{d-1} \times \mathbb{R}_+$, and the above definition will be used for the case $X = \mathbb{S}^{d-1}$, $\rho = \nu_1$ (see Section 5).

**Theorem 6.2** $(d > 1)$. *Assume that $\eta \in M^1(H)$ satisfies hypothesis $(C)$ and condition i.p. above. Then for any Borel function $\varphi$ on $V \setminus \{0\}$ with $\|v\|^{-\chi}\varphi(v)$ $\nu_1$-directly Riemann integrable, we have*

$$\lim_{t \to +\infty} t^\chi(\delta_{t^{-1}} * \nu)(\varphi) = (\nu_\infty \otimes \ell^\chi)(\varphi),$$

*where $\nu_\infty$ is a non-zero measure on $\mathbb{S}^{d-1}$. There are only 3 cases*

<u>*Case I*</u>: *$T_\mu$ has no convex invariant cone,*
   *then $Co\ (supp\ \nu) = \mathbb{R}^d$ and there exists $C > 0$ such that $\nu_\infty = C\nu_1$.*

<u>*Case II*</u> : *$T_\mu$ has a convex invariant cone $\mathcal{C} \subset \mathbb{R}^d$.*

<u>*Case II 1*</u>: *$\overline{\Lambda^a(T_\eta)} \cap \mathbb{S}^{d-1}_\infty \supset \Lambda_1(T)$,*
   *then $Co(supp\ \nu) = \mathbb{R}^d$ and there exist $C_+ > 0$, $C_- > 0$ with $\nu_\infty = C_+\nu_1^+ + C_-\nu_1^-$.*

<u>*Case II 2*</u>: *$\overline{\Lambda^a(T_\eta)} \cap \mathbb{S}^{d-1}_\infty \not\supset \Lambda_1^-(T)$,*
   *then $Co(supp\ \nu) \neq \mathbb{R}^d$ and there exists $C_+ > 0$ with $\nu_\infty = C_+\nu_1^+$.*

For a subset $Y \subset \mathbb{S}^{d-1}$, we define the polar subset $Y^\perp$ by

$$Y^\perp = \{u \in \mathbb{S}^{d-1} : \forall x \in Y, < u, x > \leq 0\}.$$

Then $Co(Y^\perp) = Y^\perp = (CoY)^\perp$ and $Y^\perp = \phi$ if $Y$ is symmetric.

With these notation we have

**Corollary 6.3.** *With the above hypothesis, there exists a continuous function $c(u)$ on $\mathbb{S}^{d-1}$, such that*

$$\forall\ u \notin (supp\ \nu_\infty)^\perp, \qquad \lim_{t \to +\infty} t^\chi \nu(H_t^u) = c(u) > 0 \ .$$

*Furthermore the set of zeros of $c$ is equal to $(supp\ \nu_\infty)^\perp$.*

In particular, the conditions $Co(supp\ \nu_\infty) = \mathbb{S}^{d-1}_\infty$ and $Co(supp\ \nu) = \mathbb{R}^d$ are equivalent.

This corollary extends the results of [9], [16], [17] to the general case.

A simplified situation is described in the following (see [6]).

**Corollary 6.4.** *Assume that for any $s > 0$, $\int \|a\|^s d\mu(a) < +\infty$ and $\int \|b\|^s d\overline{\eta}(b) < +\infty$. Also assume that $\alpha = \lim_n \frac{1}{n} \int \log \|a\| d\mu^n(a) < 0$, $T_\eta$ has no fixed point in $V$, $T_\mu$ contains a quasi-expanding element, $T_\mu$ is Zariski dense in $G$ and do not preserve a convex cone in $\mathbb{R}^d$.*

*Then, there exist $\chi > 0$ and $c(u) > 0$, such that*

$$\forall\ u \in \mathbb{S}^{d-1}, \lim_{t \to +\infty} t^\chi \nu(H_t^u) = c(u) \ .$$

*Furthermore, $Co(supp\ \nu) = \mathbb{R}^d$.*

In particular, if $\eta$ has a density on $H$ which is non-zero at $e = (0, Id) \in H$, $\alpha < 0$, and for any $s > 0$, $\int [\|a\|^s + \|b\|^s] d\eta(h) < +\infty$, then the conditions of corollary 6.4 are satisfied, hence $c(u) > 0$ on $\mathbb{S}^{d-1}$ and $supp\ \nu_\infty = \mathbb{S}^{d-1}$. Situations of this type were considered in [16], [17].

### 6.2. Remarks

(a) It can be seen that hypothesis $(C)$ in the theorem is necessary for the validity of the first conclusion.

Condition i.p. of Definition 5.3 is not necessary, but our set of conditions is generically satisfied, as we explain now. If $\eta$ is of the form $\eta = \sum_1^r p_i \delta_{h_i}$ with $r \geq 2$, $p_i > 0$, $h_i \in H$, and $h^{(r)} = (h_1, h_2, \cdots h_r)$ varies in a certain Zariski open subset $U$ of $H^{(r)}$, then all the conditions of the theorem are satisfied. Hence our conditions are generically satisfied in the weak topology of measures on $H$, and stability of the conclusions under perturbation is valid. In particular, and in contrast to the case $d = 1$, Diophantine conditions do not appear explicitly in Theorem 6.2 and corollaries 6.2 and 6.4. This is a consequence of Proposition 5.7.

(b) However, interesting special cases are not covered by the theorem. For example, if $T_\mu$ consists of diagonal matrices with real or complex entries, then condition i.p. is violated. Estimation at infinity of Poisson kernels on homogeneous nilpotent groups also leads to such problems and to analogous results (See [4]).

(c) In the simple example mentioned in section 2 $(d = 2)$, i.e., $\eta = p\delta_h + q'\delta_{h'}$ with $p, p' > 0$, and

$$h = \rho \begin{pmatrix} \cos\theta & -\sin\theta \\ \sin\theta & \cos\theta \end{pmatrix}, \quad h' = \left[ b, \begin{pmatrix} \lambda & 0 \\ 0 & \lambda' \end{pmatrix} \right],$$

where $b \neq 0$, $\rho < 1$, $\theta/\pi \notin \mathbb{Q}$ and $\lambda > 1 > \lambda' > 0$, the condition $\alpha < 0$ is satisfied for $\rho$ sufficiently small. Then $\lim_{s \to +\infty} k(s) = +\infty$, the conditions of the theorem are satisfied and we are in case I. In particular, the support of $\nu_\infty$ is the whole unit circle, and $c(u)$ is positive for every $u$. Also $Co(supp\ \nu) = \mathbb{R}^2$.

(d) In general $\overline{supp\ \nu} \cap \mathbb{S}^{d-1}_\infty$ is strictly larger than $\overline{supp\ (\nu_\infty \otimes \ell^\chi)} \cap \mathbb{S}^{d-1}_\infty = supp\ \nu_\infty$. The set $supp\ (\nu_\infty \otimes \ell^\chi)$ can be thought of as a kind of "spine at infinity" of $supp\ \nu$. In general, and in contrast to $supp\ \nu$, it has a transversal Cantor structure given by $supp\ \nu_\infty$ (see [13,14]).

(e) We have assumed $\det(a) \neq 0$ in order to rely on the group framework of [11], [13]. However, the above statements remain valid if $H$ is replaced by the affine semi-group of $V$, i.e., if $\det(a) = 0$ is allowed, as in Theorem 3.1.

## 7. Some tools of the proof

### 7.1. The scheme

We consider the functional equation

$$(1) \qquad\qquad z - b = a\ (z \circ \theta)\ .$$

Let $\nu'$ he the law of $z - b$. Then we obtain

$$\nu' = \mu * \nu, \quad \nu - \nu' = \nu - \mu * \nu\ .$$

It can be shown that $\nu$ is given by the potential

$$\nu = \sum_0^\infty \mu^k * (\nu - \nu')\ .$$

The advantage of this implicit formula is that it involves only convolution of measures on $GL(V)$ and $V \setminus \{0\}$. It can be shown that $\nu - \nu'$ is "small at infinity", then the general renewal theorem of [15] can be applied to the operator $P_\mu$ on $V \setminus \{0\} = \mathbb{S}^{d-1} \times \mathbb{R}_+$ defined by

$$P_\mu(v, .) = \mu * \delta_v \ ,$$

and to its potential kernel

$$\sum_0^\infty P_\mu^k(v, .) = \sum_0^\infty \mu^k * \delta_v \ .$$

Thus the invariant measures $\lambda$ of $P_\mu$, defined by $\lambda = \mu * \lambda$, play an essential role in the problem (see Section 4). A technically important step is to consider the space $\tilde{V} = V/\{\pm Id\}$, since the geometric properties of convolution operators on $\tilde{V}$ are simpler than on $V$. In order to verify the conditions of validity for the renewal theorem of [15], we need to use the following properties, which are developed in [11], [13].

(a) The equation

$$r^\chi(\mu)(\rho) = \rho \ , \qquad \rho \in M^1(\mathbb{P}^{d-1}),$$

has a unique solution on $\mathbb{P}^{d-1}$.

(b) The spectrum

$$\Sigma(T) = \{\log[r(a)] : a \in T_\mu^{prox}\}$$

generates a dense subgroup of $\mathbb{R}$ (see also [14]).

(c) The Markov chain on $\mathbb{P}^{d-1}$ associated with the above equation has strong mixing properties and the first Lyapunov exponent of $\mu$ is simple.

These properties depend on the condition i.p., of Definition 5.3 which justifies its introduction.

(d) Finally we need to study the positivity of $C, C_+$ and $C_-$.

Here we use the geometry of $T_\eta$ and $T_\mu$, and we study an auxiliary Markov chain on the space of affine hyperplanes of $\mathbb{R}^d$. This allows us to control the function $\nu(H_t^u)$ using the ergodic properties of this chain (see [18]).

### 7.2. *Analytic argument for $d = 1$*

For $d = 1$, a special proof of the positivity of $c_+$ and $c_-$ can be given under some stronger hypothesis than in Theorem 3.1. We sketch it here for case II under the simplifying condition:

$$|b| \leq B, \quad s_\infty = \infty, \quad a > 0 \ .$$

We restrict to the study of $c_+$. Then, the above functional equation gives

$$(z - b)^+ = a(z^+ \circ \theta) \ ,$$

where $z^+ = \sup\{z, 0\}$. We denote by $\nu_+$ the law of $z^+$, and write $h(s) = \mathbb{E}(|z^+|^s)$, $u(s) = \mathbb{E}(|z^+|^s - |(z - b)^+|^s)$. Observe that $h$ is defined for $s < \chi$ while $|u(s)| \leq B^s$ for any $s$, hence the function $\frac{u(s)}{1 - k(s)}$ is meromorphic in the half plane $\operatorname{Re}(s) > 0$.

The above equation gives for $0 \leq s < \chi$

$$h(s)[1 - k(s)] = u(s) \ ,$$

and from the renewal theorem (see [9], [10]),

$$\lim_{t \to +\infty} t^\chi (\delta_{t^{-1}} * \nu_+) = C_+ \ell^\chi \ ,$$

where $C_+$ is the residue at $\chi$ of $\frac{u(s)}{1-k(s)}$ (see [5]).

If $C_+ = 0$, then it follows that the function $\frac{u(s)}{1-k(s)}$ is holomorphic at $\chi$, and hence on the whole line $\mathbb{R}_+$.

We recall that the Mellin transform $\widehat{\gamma}(s) = \int_0^\infty x^s d\gamma(x)$ of a positive measure on $\mathbb{R}_+$ cannot be extended holomorphically to a neighbourhood of its abcisse of convergence $\tau$ ([19] p. 58). In particular, this is valid for $u(s) = \mathbb{E}(|z^+|^s)$, hence the condition $C_+ = 0$ implies that $E(|z^+|^s) < +\infty$ in a neighborhood of $\chi$. The same argument gives

$$\tau = +\infty, \quad [k(s) - 1]E(|z^+|^s) \le B^s \quad (s > \chi).$$

It follows that

$$|a|_\infty |z^+|_\infty = \lim_{s \to +\infty}[k(s)]^{1/s}\mathbb{E}(|z^+|^s)^{1/s} \le B \ .$$

Since $|a|_\infty > 0$, this implies that $+\infty \notin \overline{\Lambda^a(T)}$, which contradicts assumptions II 1, II 2.